\documentclass[a4paper,12pt,final]{amsart}
\usepackage{times,a4wide,mathrsfs,amssymb,amsmath,amsthm,enumerate,xypic,tikzsymbols,dsfont}

\newcommand{\C}{\mathbb{C}}
\newcommand{\ZZ}{\mathbb{Z}}

\newcommand{\QQ}{\mathbb{Q}}
\newcommand{\NN}{\mathbb{N}}
\newcommand{\PP}{\mathbb{P}}

\newcommand{\OO}{\mathcal O}
\newcommand{\Ss}{\mathcal S}

\newcommand{\XX}{\mathcal X}
\newcommand{\YY}{\mathcal Y}

\newcommand{\pic}{\hbox{Pic}}

\newcommand{\rom}{\romannumeral}

\DeclareMathOperator{\aut}{Aut}

\DeclareMathOperator{\ide}{id}

\DeclareMathOperator{\ima}{Im}

\DeclareMathOperator{\GDCH}{GDCH}
\DeclareMathOperator{\CH}{CH}
\newtheorem{theorem}{Theorem}[section]

\newtheorem{lemma}[theorem]{Lemma}

\newtheorem{proposition}[theorem]{Proposition}
\newtheorem{conjecture}[theorem]{Conjecture}
\newtheorem{remark}[theorem]{Remark}
\newtheorem{definition}[theorem]{Definition}
\newtheorem{convention}{Conventions}

\newtheorem{nonumbering}{Theorem}

\newtheorem{nonumberingc}{Conjecture}

\newtheorem{nonumberingt}{Acknowledgements}

\newtheorem{nonumberingcon}{Statements and declarations}

\begin{document}

\author[Gilberto Bini]
{Gilberto Bini}

\address{Dipartimento di Matematica e Informatica, \
Via Archirafi 34 -- 90138 Palermo, \
Universit\`a degli Studi di Palermo, }
\email{gilberto.bini@unipa.it}

\author[Robert Laterveer]
{Robert Laterveer}

\address{Institut de Recherche Math\'ematique Avanc\'ee,
CNRS -- Universit\'e 
de Strasbourg,\
7 Rue Ren\'e Des\-car\-tes, 67084 Strasbourg CEDEX,
FRANCE.}
\email{robert.laterveer@math.unistra.fr}

\title{Zero-cycles and the Cayley--Oguiso automorphism}

\begin{abstract} Cayley and Oguiso have constructed certain quartic K3 surfaces $S$, with automorphisms $g$ of infinite order. We show that when $g$ is symplectic (resp. anti-symplectic), it acts as the identity (resp. minus the identity) on the degree zero part of the Chow group of zero-cycles of $S$. 
\end{abstract}

\thanks{GB and RL both gratefully acknowledge the support of Laboratoire Ypatia des Sciences Mathématiques LYSM CNRS-INdAM International Research Laboratory.
RL is also supported by the Agence Nationale de Recherche under grant ANR-20-CE40-0023.}

\maketitle

\section{Introduction}

Given a smooth projective variety $Y$ over $\C$, let $\CH^i(Y)_{\ZZ}$ denote the Chow groups of $Y$ (i.e. the groups of codimension $i$ algebraic cycles on $Y$ modulo rational equivalence), and let
$\CH^i_{hom}(Y)_{\ZZ}$ denote the subgroup of homologically trivial cycles.

We are interested in the following conjecture, which is a relative version of the famous Bloch conjecture (see, for instance, \cite{bloch}). 

\begin{nonumberingc}[Conjecture \ref{conj}] Let $S$ be a K3 surface, and let $g\in\aut(S)$. 
 
 \noindent
 (\rom1) If $g$ is symplectic then $g$ acts as the identity on $\CH^2_{hom}(S)_{\ZZ}$.
 
 \noindent
 (\rom2) If $g$ is anti-symplectic then $g$ acts as minus the identity on $\CH^2_{hom}(S)_{\ZZ}$.
 \end{nonumberingc}

For the state of the art on this conjecture, we refer to Section \ref{conjecture}. After seminal work by Huybrechts and Voisin, the conjecture remains open for K3 surfaces with Picard number $2$ and automorphism group isomorphic to the additive group of integers. Our main result settles Conjecture \ref{conj} for an 18-dimensional family of K3 surfaces with infinite-order automorphisms. More precisely, the following holds.

\begin{nonumbering}[Theorem \ref{main}] Let $S$ be a Cayley--Oguiso K3 surface, and let $g \in\aut(S)\cong\ZZ$ be any automorphism. Then Conjecture \ref{conj} is true for $(S,g)$.
\end{nonumbering}

The Cayley--Oguiso K3 surfaces were first introduced by Cayley in \cite{Cay} as smooth determinantal quartics in three-dimensional projective space. More precisely, given $64$ complex numbers $a_{ijk}$, namely a general $4 \times 4 \times 4$ tensor in modern terms, Cayley defined smooth quartics as the zero locus of the determinant of $16$ linear forms in $4$ variables. In 2015, Oguiso described in \cite{Og1} an $18$-dimensional family of algebraic K3 surfaces with a fixed point-free automorphism $g$ of positive entropy. As shown in \cite{F+}, every K3 surface $S$ of this family can be embedded in three-dimensional projective space, thus yielding an infinite sequence $S_n$ of K3 surfaces isomorphic to one another as abstract manifolds. What's more, the automorphism group of $S$ is an infinite cyclic group generated by $g$, which is anti-symplectic. Every surface $S_n$ carries a smooth curve of genus $3$ and degree $6$, so $S_n$ is a determinantal quartic, as proved in \cite{Beau1}. As such, they can be described as Cayley did in 1870!

More interestingly, $S_n$ and $S_{n+3r}$ are not only isomorphic as abstract K3 surfaces, but also as quartics in three dimensional projective space via a projective realization of $g$. On the other hand, for every integer $n$ the quartics $S_n$, $S_{n+r}$, $S_{n+2r}$ are not pairwise projectively equivalent, as shown in \cite{Og5}, where $S$ is embedded in a product of two three dimensional projective spaces as a complete intersection of four hypersurfaces of type $(1,1)$. Anyhow, as proved in loc. cit., $S_n$, $S_{n+r}$, $S_{n+2r}$ are pairwise isomorphic via a Cremona transformation, i.e., they are Cremona isomorphic. Therefore, $S$ produces three related quartic surfaces. This is somewhat reminiscent of the $SO(8)$-triality phenomenon (see Subsection \ref{triality} for more details), which is further investigated in \cite{manivel} from an algebraic geometry viewpoint and more generally applied to hyperk\"ahler varieties.

In Section \ref{cayleyoguiso}, we briefly survey the construction, as well as the properties, of the Cayley-Oguiso K3 surface $S$. In Section \ref{s: auto}, we list some properties of the automorphism $g$ on $S$. In particular, we recall that $g$ induces an automorphism $g_{ind}$ on the Hilbert scheme $S^{[2]}$ of $0$-dimensional subscheme of length $2$. As first pointed out by Oguiso in \cite{Og2}, $g_{ind}$ is equal to the composition $\iota_0 \circ \iota_1 \circ \iota_2$, where $\iota_l$ is the Beauville involution on $S_l^{[2]} \simeq S^{[2]}$ for $l=0,1,2$. In Section \ref{otherquartics}, we recall the K3 surfaces $S_{m,a}$ introduced in \cite{Lee} where $m$ is an integer greater than or equal to $2$ and $a$ is a positive integer. For $(m,a)=(2,1)$ the surface $S_{2,1}$ is the Cayley-Oguiso surface. In general, $S_{2,a}$ is a quartic K3 surface for any $a$ and it is determinantal only for $a=1$. Moreover, for any pair $(m,a)$ the automorphism group is infinite cyclic. In Section \ref{conjecture}, we give a brief overview on the state of the art on Conjecture \ref{conj}. In Section \ref{proofs}, we include two different proofs of Theorem \ref{main}. The first proof involves the Franchetta property, which is recalled in Section \ref{franchetta1}. The second proof uses the relation with the Beauville involution on the Hilbert square, and a result on the action of this involution on zero-cycles. For these purposes, we recall the structure of the Chow ring of a K3 surface in Section \ref{s:bv}. Finally, in Section \ref{otherfamilies}, we also settle the conjecture for two other 18-dimensional families of K3 surfaces, that is, $S_{2,2}$ and $S_{2,3}$: see Theorem \ref{main2}. We hope to investigate further the surfaces $S_{2,a}$ for $a \geq 4$ for new examples of K3 surfaces where Conjecture \ref{conj} might hold.

 \vskip0.6cm

\begin{convention} In this paper, the word {\sl ``variety''\/} will refer to a reduced irreducible scheme of finite type over $\C$. A {\sl subvariety\/} is a (possibly reducible) reduced subscheme which is equidimensional. 

{\bf All Chow groups will be with rational coefficients}: we will denote by $\CH_j(Y)$ the Chow group of $j$-dimensional cycles on $Y$ with $\QQ$-coefficients; for $Y$ smooth of dimension $n$ the notations $\CH_j(Y)$ and $\CH^{n-j}(Y)$ are used interchangeably. 
The notation $\CH^j_{hom}(Y)$ will be used to indicate the subgroup of homologically trivial cycles.
For a morphism $f\colon X\to Y$, we will write $\Gamma_f\in A_\ast(X\times Y)$ for the graph of $f$.

\end{convention}

\section{The Cayley--Oguiso Surfaces}\label{cayleyoguiso}

\subsection{Quartic Surfaces in ${\mathbb P}^3$}

Let $X$ and $Y$ be two degree $d$ smooth hypersurfaces in ${\mathbb P}^n$ that are isomorphic as abstract varieties. It is then natural to ask if $X$ and $Y$ are {\em projectively equivalent} or {\em Cremona equivalent}. In the former case, there exists an automorphism of ${\mathbb P}^n$ that restricts to an isomorphism of $X$ and $Y$. In the latter case, there exists a Cremona transformation $f: {\mathbb P} \dashrightarrow {\mathbb P}$ such that the restriction to $X$ extends to an isomorphism onto $Y$. 

If $(n,d) \neq (3,4), (2,3)$, Matsumura and Monsky proved that $X$ and $Y$ are projectively equivalent (so Cremona equivalent) if they are isomorphic as abstract varieties \cite{mm64}. As a consequence, the case of smooth quartic surfaces in $3$-dimensional projective space is particularly relevant. 

For these purposes, we focus on a construction due to K. Oguiso in \cite{Og2} (cfr. \cite{F+}), who proves that there exists a dense subset of an $18$-dimensional family of K3 surfaces $S$ with N\'eron-Severi group isomorphic to $(N,b)$. Here $N$ is the ring of integers of ${\mathbb Q}(\eta)$ for $\eta^2=1+\eta$ and $b$ is the bilinear form given by
$$
\left( 
\begin{array}{cc}
4 & 2 \\
2 & -4
\end{array}
\right)
$$

In \cite[Theorem 1.15]{Og5}, K. Oguiso proves that $S$ has two very ample divisors $h_1$ and $h_2$. The map induced by the linear series $|h_1|$ gives a smooth quartic surface $S_1$ in ${\mathbb P}^3$ that is determinantal, i.e., it is the zero locus of the determinant of a matrix $T_1(x)$ with entries given by $16$ linear forms in four variables. The very ample divisor $h_2$ gives a genus $3$ curve of degree $6$; hence $S$ is embedded in ${\mathbb P}^3$ as a determinantal quartic: see \cite{Beau1}.

As proved in \cite[Theorem 1.15]{Og5}, it turns out that there is an isomorphism
\[S \simeq Q_1 \cap Q_2 \cap Q_3 \cap Q_4 \subset P\ ,\] 
where $Q_j$ are hypersurfaces of bidegree $(1,1)$ of $P:={\mathbb P}^3 \times {\mathbb P}^3$. Moreover, if we set $V= Q_1 \cap Q_2 \cap Q_3 \subset P$ for any $S$, the birational map given by the composition ${p_2}_{|V} \circ ({p_1}_{|V})^{-1}$ shows that $S_1$ and $S_2$ are Cremona equivalent, where $p_i$ is the projection onto the $i$-th factor of $P$. What's more, for very general hypersurfaces $Q_j$ the surfaces $S_1$ and $S_2$ are Cremona isomorphic but not projectively equivalent. As proved by F. Reede in \cite{reede}, the birational isomorphism is actually given by the cubo-cubic transformation in ${\mathbb P}^3$.

\subsection{Quartic Determinantal Surfaces from a Tritensor}
\label{tritensor}

In the previous section we started from a smooth K3 surface $S$ and defined two smooth quartic surfaces $S_1$ and $S_2$ in ${\mathbb P}^3$, both of which are determinantal quartics. As observed by Cayley \cite{Cay} and further developed by Festi, Garbagnati, van Geemen and van Luijk \cite{F+}, given a determinantal quartic, it is possible to find two others, not only one.

Let $W$ be a complex vector space of dimension $4$. Let $e_0, \ldots, e_3$ be a basis of $V$. Pick an element $T \in W \otimes W \otimes W$ and write it as $w= \sum_{i,j,k=0}^3 m_{i,j,k}e_i \otimes e_j \otimes e_k$ where $m_{i,j,k}$ are complex numbers. Fix anyone of the indices $i,j,k$ and define the $4 \times 16$ {\em flattening matrices} $F_0$, $F_1$, $F_2$ of $T$ as follows:
$$
F_0=[m_{0,j,k}|m_{1,j,k}|m_{2,j,k}| m_{3,j,k}],
$$
$$ 
F_1=[m_{i,0,k}|m_{i,1,k}|m_{i,2,k}| m_{i,3,k}]
$$
$$
F_2=[m_{i,j,0}|m_{i,j,1}|m_{i,j,2}| m_{i,j,3}].
$$

Fix a basis of $V$ and denote by $x_0, x_1, x_2, x_3$ the homogeneous coordinates on ${\mathbb P}(V)$. If we multiply each of the flattening matrices by the $16 \times 4$ matrix 
$$
\left[
\begin{array}{c}
x_0 I_4 \\
x_1 I_4 \\
x_2 I_4 \\
x_3 I_4
\end{array}
\right]
$$
we obtain the three $4 \times 4$ matrices 
$$
T_0(\underline{x})=\left(  \sum_{i=0}^3 x_i m_{i, j, k}\right)_{j,k} , \quad T_1(\underline{x})=\left( \sum_{j=0}^3 x_j m_{i,j,k} \right)_{i,k}, \quad T_2(\underline{x})= \left( \sum_{k=0}^3 x_k m_{i,j,k} \right)_{i,j}.
$$
Therefore, from any tritensor $T$ we can define three quartic surfaces $\Sigma_l$ in ${\mathbb P}(V)$ as $\det(T_l(\underline{x}))=0$ for $l=0,1,2$. These surfaces are related as follows. For example, start with $\Sigma_0$. For every point $[x_0: x_1: x_2: x_3] \in \Sigma_1$ there exists a non-zero solution (which is unique up to non-zero multiples) $Y(\underline{x})$ of the linear system $T_1(\underline{x})Y(\underline{x})=0$. Denote by $y_k(\underline{x})$ the components of $Y(\underline{x})$ for $k=0,1,2,3$. Therefore there exists a well-defined morphism $\varphi_{12}: \Sigma_0 \to {\mathbb P}(V)$ such that 
$
\varphi_{12}([x_0: x_1: x_2: x_3]) =[y_0(\underline{x}): y_1(\underline{x}): y_2(\underline{x}): y_3(\underline{x})].
$
The image $\varphi_{01}(\Sigma_0)$ is equal to $\Sigma_1$. For these purposes, recall that Laplace's rule implies that for $k=0,1,2,3$ the component $y_k(\underline{x})$ is given by $(-1)^kT_{0,k}(\underline{x})$, i.e. the algebraic complements of the elements of the first row of $T_0$. Therefore we have $\det(T_1(Y(\underline{x})))=0$. As A. Cayley wrote in \cite{Cay}, the process can be repeated indefinitely. In other words, we can repeat a similar construction starting from $S_1$ or from $S_2$, thus obtaining the morphisms $\varphi_{12}: \Sigma_1 \to \Sigma_2$ and $\psi: \Sigma_2 \to {\mathbb P}(V)$. By the same arguments as before, the image of $\psi$ is $\Sigma_0$ and, consistently, we will write $\psi=\varphi_{20}$.

\subsection{Triality Phenomena}\label{triality}

Let $W_1$, $W_2$, $W_3$ be vector spaces over the real field. A triality among them is a trilinear map $\tau: W_1 \times W_2 \times W_3 \to {\mathbb R}$ that is non-degenerate, i.e., if we fix any two arguments to any non-zero values, the linear functional induced on the third space is nonzero. 

Let $Spin(8)$ be the spinor group of real $8 \times 8$ matrices, which is a double cover of $SO(8)$. Most of what we are going to recall holds for $n \geq 8$ (see, for instance, \cite{baez}) but here we focus on $n=8$ as it is the case that is of our interest. The group $Spin(8)$ has a natural representation, namely a vector space $V_8$, as well as two spinor representations $S_8^+$ and $S_8^-$, both of which are $8$-dimensional. As recalled in loc. cit., it is possible to define a trilinear map $t_8: V_8 \times S_8^{+} \times S_8^{-} \to {\mathbb R}$ that is in fact a triality. Moreover,  the group $Out(Spin(8))$ of outer automorphisms of $Spin(8)$ is isomorphic to the symmetric group of order $6$. As a consequence, the three representations above are permuted one onto another by linear maps that are induced by outer automorphisms of $Spin(8)$.

Finally, recall that $Spin(8)$ has a non-faithful action on $S^7$, the unit sphere in ${\mathbb R}^8$. This action descends on ${\mathbb P}(W_i^*)$ for $W_1=V_8$, $W_2=S_8^+$ and $W_3=S_8^-$, where ${\mathbb P}(W^*_i)$ may be viewed as $3$-dimensional complex projective space. By the preceding arguments, triality permutes these three projective spaces. Since the permutations are induced by linear maps, any surface in one of these projective spaces is mapped to a surface (of the same degree) in another of these projective spaces. Therefore, we can consider the construction in Subsection \ref{tritensor} as a {\em triality phenomenon}: for further reading we refer, for instance, to \cite{manivel} and to \cite{BH} where this phenomenon is studied in depth (the Cayley--Oguiso surfaces occur in \cite{BH} under the appelation ``Rubik's revenge'', cf. \cite[Section 4]{BH}). 

\section{The Automorphism Group of the Cayley--Oguiso Surfaces}
\label{s: auto}

As proved in \cite[Theorem 4.1]{Og6}, the surface $S$ introduced before admits an automorphism $g$ without fixed points and of positive entropy. In other words, let $g^*$ be the ${\mathbb C}$-linear extension on the degree $2$ singular cohomology of $S$. Then $g$ is of positive entropy if each eigenvalue $\lambda$ of $g^*$  has absolute value greater than $1$. Moreover, as further developed in \cite[Theorem 1.2]{F+}, the automorphism group of $S$ is isomorphic to the group of integers. Any of its generators is fixed point free and has positive entropy. In particular, there exists a generator $g$ such that the action on the transcendental lattice $T(S)$ is multiplication by $-1$ and on the N\'eron-Severi group $NS(S)\simeq {\mathbb Z}[\eta]$ is multiplication by $\eta^6$ where $\eta^2=1+\eta$.

Let $D_n$ be the very ample divisors in $NS(S)$ that corresponds to $\eta^{2n}$. Every $D_n$ embeds $S$, via an embedding $\phi_n$, onto $S_n$ in ${\mathbb P}^3$ because $D_n^2=4$. The quartics $S_n$ are determinantal. Moreover, $S_{n + 3r} \simeq \Sigma_n $ for $r=0,1,2$. Indeed, 
$g^*(D_0)=D_3$ and the corresponding surfaces $S_0$ and $S_3$ are mapped one onto the other by a projective transformations of ${\mathbb P}^3$. As remarked in \cite[\S 1.10]{F+}, we have the equality \[ g=\phi_0^{-1} \circ \phi_3\ ,\] 
which implies that $S_0$ and $S_3$ are projectively equivalent, as well as isomorphic as abstract varieties. On the other hand, the automorphism $g$, and the action $g^*$ on the N\'eron-Severi group, does not induce a projective transformation between any of the pairs $S_i$, $S_j$ for $i \neq j \in \{0,1,2\}$. 

We finally remark that in loc. cit., Section 4, a concrete description of $g$ is given in terms of polynomial maps.

\subsection{Relation with Beauville's involution}\label{beauvilleinvolution}

For any quartic K3 surface $S$, Beauville has constructed an involution $\iota$ of the Hilbert square $S^{[2]}$ \cite{Beau}. By definition, the involution $\iota$ sends an unordered pair of points $(x_1,x_2)$ (where $x_i\in S$) to the residual intersection of the line spanned by $x_1,x_2$ and the quartic $S$. The involution $\iota$ is a birational map in general; it is an automorphism when $S$ does not contain lines.

\begin{proposition}[Oguiso \cite{Og2}]\label{og} Let $S$ be a Cayley--Oguiso quartic, and let $g\in\aut(S)\cong\ZZ$ be the generator as above. Then $g$ induces an automorphism $g_{ind}$ of $S^{[2]}$, and there is a relation
  \[  g_{ind}=\iota_2\circ\iota_1\circ\iota_0\ ,\]
  where $\iota_l$ indicates the Beauville involution of $S_l$ (and $S_l\subset\PP^3$ is as above embedded by the divisor $D_l$).
\end{proposition}    

\begin{proof} This is part of \cite[Example 1 in Section 5]{Og2}.

\end{proof}

  \section{Other quartics with infinite-order automorphisms}\label{otherquartics}
  
  Recently, Lee has discovered a generalization of the Cayley--Oguiso construction:
  
  \begin{theorem}[Lee \cite{Lee}]\label{lee} Let $S=S_{m,a}$ be the K3 surface with N\'eron--Severi lattice
    \[  m  \left( 
\begin{array}{cc}
2 & a \\
a & -2
\end{array}
\right)
   ,\]
    where $(m,a)\in\NN_{\ge 2}\times\NN_{\ge 1}$.
  
\noindent
(\rom1)  
  $\aut(S)\cong\ZZ$. 
  
  \noindent
  (\rom2) If $m$ divides $a$, then $\aut(S)$ is generated by 
  a symplectic automorphism $h$. If $m$ divides $a^2+1$, then $\aut(S)$ is generated by an anti-symplectic automorphism $g$.
  
  \noindent
  (\rom3) If $m=2$, $S$ is a quartic $K3$.
   \end{theorem}
   
   \begin{proof} This is \cite[Theorem 1 and Section 5.1]{Lee}.
    \end{proof}
   
   \begin{remark} The Cayley--Oguiso surfaces correspond to $(m,a)=(2,1)$. The surfaces $S_{2,a}$ with $a>1$ are not determinantal quartics, because they do not contain a sextic curve of genus 3 (cf. \cite[Remark 1]{Lee}).
    \end{remark}
   
   Similar to Proposition \ref{og}, there is a relation with Beauville involutions:
   
   \begin{proposition}[Lee \cite{Lee}]\label{lee2} Let  $S=S_{2,a}$ and $g, h$ be as in Theorem \ref{lee}. Then $g$ induces an automorphism $g_{ind}$ of $S^{[2]}$, and there is a relation
  \[   \begin{split}  &\iota_{2\ell +k} = (h_{ind})^\ell\circ \iota_k\circ  (h_{ind})^{-\ell} \ \ \ \ (0\le k\le 1)\ \ \ \ \hbox{if\ $a$\ is\ even}\ ;\\
                                 &\iota_{3\ell +k} = (g_{ind})^\ell\circ \iota_k\circ  (g_{ind})^{-\ell} \ \ \ \ (0\le k\le 2)\ \ \ \ \hbox{if\ $a$\ is\ odd}\ .\\  
                                 \end{split}\]
                                 \end{proposition}
                                 
                                 \begin{proof} This is \cite[Theorem 2]{Lee}.
                                 \end{proof}

 \section{The conjecture}\label{conjecture}
 
 The following conjecture can be seen as a relative version of the famous Bloch conjecture: 
 
 \begin{conjecture}\label{conj} Let $S$ be a K3 surface, and let $g\in\aut(S)$. 
 
 \noindent
 (\rom1) If $g$ is symplectic then $g$ acts as the identity on $\CH^2_{hom}(S)_{\ZZ}$.
 
 \noindent
 (\rom2) If $g$ is anti-symplectic then $g$ acts as minus the identity on $\CH^2_{hom}(S)_{\ZZ}$.
 \end{conjecture}
 
 \begin{remark} It follows from Roitman's theorem \cite{Roit} that $\CH^2_{hom}(S)_{\ZZ}$ is torsion-free for a K3 surface $S$, and so the conjecture is equivalent to the version where $\CH^2_{hom}(S)_{\ZZ}$ is replaced by the $\QQ$-vector space $\CH^2_{hom}(S)_{}$.
   \end{remark}

 \begin{theorem}[Huybrechts \cite{Huy}]\label{huy} Conjecture \ref{conj} is true when $g$ is a finite-order symplectic automorphism.
 \end{theorem}
 
 \begin{proof} This is \cite[Theorem 0.1]{Huy}, which uses derived category techniques and builds on work of Voisin concerning involutions \cite{Vo1}.
 
 \end{proof}

%
%

  \begin{theorem}[Du--Liu \cite{DL}]\label{dl} Conjecture \ref{conj} is true when $S$ is an elliptic K3 surface, and $g\in\aut(S)$ is a symplectic automorphism preserving the elliptic fibration.
   \end{theorem}
   
   \begin{proof} This is \cite[Theorem 1.5]{DL}. 
       \end{proof}
 
 \begin{remark} Examples of $S$ and $g$ as in Theorem \ref{dl} with $g$ of infinite order are contained in \cite{Nik1}, \cite{Nik2}.
 
 \end{remark}

 \begin{theorem}[Pawar \cite{Paw}] Conjecture \ref{conj} is true when $S$ is a Kummer K3 surface, and $g\in\aut(S)$ is a symplectic automorphism induced by an automorphism from the abelian surface. 
 \end{theorem}
 
 \begin{proof} This is \cite[Theorem 5.2.2]{Paw}.
  \end{proof}
  
  \begin{theorem}[Li--Yu--Zhang \cite{LYZ}]\label{lyz} Conjecture \ref{conj} is true when $S$ has Picard number at least 3.
  \end{theorem}
  
  \begin{proof} This is \cite[Theorem 1.5]{LYZ}. The approach is via derived categories; the paper considers more generally the action of auto-equivalences on the Chow group of zero-cycles.
  \end{proof}

  \begin{proposition}\label{easy} Conjecture \ref{conj} is true when $\aut(S)=\ZZ_2\ast\ZZ_2$.
\end{proposition}

\begin{proof} Let $\sigma_1, \sigma_2$ denote the two generators. Conjecture \ref{conj} is true for $\sigma_1$ and for $\sigma_2$ thanks to Theorem \ref{huy}. But then the conjecture is also true for any composition of $\sigma_1$ and $\sigma_2$.
\end{proof}

\begin{remark} For K3 surfaces of Picard number 2, it is known that the only infinite automorphism groups that can occur are $\ZZ$ and $\ZZ_2\ast\ZZ_2$ \cite{GLP}.
In view of Theorem \ref{lyz} and Proposition \ref{easy}, it follows that the only open cases of Conjecture \ref{conj} are K3 surfaces $S$ of Picard number 2 and $\aut(S)=\ZZ$. Actually, as proved in \cite[Theorem 5.3]{hkl}, any fixed point free automorphism of a projective K3 surface with Picard number $2$ is the Cayley--Oguiso automorphism $g$ on $S$.
\end{remark}

  \section{The Chow ring of K3 surfaces}
  \label{s:bv}
   We recall the following famous result on the Chow ring of K3 surfaces.
  
  \begin{theorem}[Beauville--Voisin \cite{BV}] Let $S$ be a projective K3 surface. There exists a class $o_S\in \CH^2(S)_{\ZZ}$, with the following properties:
  
  \begin{enumerate}
  
  \item $D\cdot D^\prime\in \ZZ [o_S]$, for any divisors $D,D^\prime\in\pic(S)$;
  
  \item $c_2(T_S)=24 o_S$ in $\CH^2(S)_{\ZZ}$;
  
  \item $\Delta_S\cdot (p_i)^\ast(D) = D\times o_S + o_S\times D$ in $\CH^3(S\times S)_{\ZZ}$, for any $D\in\pic(S)$, where $\Delta_S\subset S\times S$ denotes the diagonal and $p_i\colon S\times S\to S$ are the two projections on the factors;
  
  \item $\Delta_S\cdot \Delta_S= 24 o_S\times o_S$ in $\CH^4(S\times S)_{\ZZ}$.
  \end{enumerate}
   \end{theorem}

 \section{Franchetta property}\label{franchetta1}

An important tool we are going to use in this paper is the Franchetta property. 

\begin{definition}\label{frank} Let $\pi\colon\YY\to B$ be a smooth projective morphism, where $\YY, B$ are smooth quasi-projective varieties. For any $b\in B$, we write $Y_b$ for the fiber $\pi^{-1}(b)$.
One says that $\YY\to B$ has the {\em Franchetta property in codimension $j$\/} if the following property holds:

\smallskip

 for every $\Gamma\in \CH^j(\YY)$ such that the restriction $\Gamma\vert_{Y_b}$ to the fiber is homologically trivial for all $b\in B$, the restriction $\Gamma\vert_{Y_b}$ is zero in $\CH^j(Y_b)$ for all $b\in B$.
 
 We will say that $\YY\to B$ has the {\em Franchetta property\/} if $\YY\to B$ has the Franchetta property in codimension $j$ for all $j$.
 \end{definition}
 
 This property has been studied in \cite{BL}, \cite{FLV}, \cite{FLV2}, \cite{FLV3}.
 
 \begin{definition}\label{def:gd} Given a family $\YY\to B$ with the properties described above, and given $Y:=Y_b$ a fiber, we write
   \[ \GDCH^j_B(Y):=\ima\Bigl( 
  \CH^j(\YY)\to \CH^j(Y)\Bigr) \]
   for the subgroup of {\em generically defined cycles}. 
  Whenever it is clear to which family we are referring, we will often suppress the index $B$ from the notation.
  \end{definition}
  
 We remark that, with this notation, the Franchetta property is equivalent to saying that $\GDCH^\ast_B(Y)$ injects into cohomology, under the cycle class map.
 
\smallskip 
 
The Franchetta property for a family $\YY \to B$ does not imply and is not implied by the Franchetta property for a subfamily $\YY^\prime \to B^\prime$, where $B^\prime$ is a closed subscheme of $B$. If $B^\prime \to B$ is a dominant morphism, the Franchetta property for the base-changed family $\XX_{B^\prime} \to B^\prime$ implies the property for $\XX \to B$ \cite[Remark 2.6]{FLV}.

 \section{Main result}\label{proofs}
 
 \begin{theorem}\label{main} Let $S$ be a Cayley--Oguiso K3 surface, and let $g$ be a generator of $\aut(S)\cong\ZZ$. Then ($g$ is anti-symplectic and) Conjecture \ref{conj} is true for $S$, i.e.
   \[ g_\ast=-\ide\colon\ \ \CH^2_{hom}(S)\ \to\ \CH^2_{hom}(S)\ .\]
   \end{theorem}
   
 \begin{proof} Without loss of generality, we may assume $g$ is the generator described in Section \ref{s: auto}.
Set $P={\mathbb P}^3 \times {\mathbb P}^3$ and introduce $V:=H^0(P, {\mathcal O}_P (1,1)^{\oplus 4})$, where ${\mathcal O}_P (1,1)$ is the tautological line bundle on $P$ of bidegree $(1,1)$.
. Also, set $\overline{B}={\mathbb P}(V)$. Let $B$ be the open set in ${\overline B}$ which parametrizes smooth complete intersection in $P$ of $4$ sections of bidegree $(1,1)$. Let $\overline{{\mathcal S}}$ be the universal complete intersection family over $B$. As such, there exist two morphisms, namely ${\overline {\mathcal S}} \times_{\overline{B}} {\overline {\mathcal S}} \to {\overline B}$ and $ \pi: {\overline {\mathcal S}} \times_{\overline{B}} {\overline {\mathcal S}} \to P \times P$. Denote by $p_i$ the projection of $P$ onto the $i$-th factor ${\mathbb P}^3$. Let $h_j\in \CH^1(S)$ denote the restriction of  $(p_j)^\ast(H)$ to $S$, where is the hyperplane divisor in ${\mathbb P}^3$.

Notice that $\pi$ is a stratified projective bundle. Indeed, for any $p \in P \simeq \Delta_P \subset P \times P$ the fiber over $p$ is a projective space of constant dimension $r$, which is related to the condition imposed by $p$ on the sections of the very ample line bundle  ${\mathcal O}_P (1,1)$. For $(p, p') \in (P \times P) \setminus \Delta_P$, the fiber over $(p,p')$ is a projective space of constant dimension $t \neq r$, as two distinct points impose two independent conditions. This stratified bundle restricts to a smooth projective morphism from the irreducible open subset ${\mathcal S} \times_B {\mathcal S}$ to $B$, the open subset of $\overline{B}$ that parametrizes smooth complete intersections.  Our next step is to prove that the Franchetta property holds for this family.

First, as proved in \cite{FLV} the stratified projective bundle argument implies that 
\[  \begin{split}
\GDCH^*_B(S_b \times S_b) &=\bigl\langle\ima\bigl(\CH^\ast(P\times P)\to \CH^\ast(S_b\times S_b)\bigr), \Delta_{S_b}\bigr\rangle\\
                       &= \bigl\langle (p_i)^*(h_j), \Delta_{S_b}\bigr\rangle\ ,\\
                          \end{split} \]
where the $\langle\cdot \rangle$ denotes the ${\mathbb Q}$-algebra generated by $p_i^*(h_j)$ and the diagonal $\Delta_{S_b}$. Second, apply the results in \cite{BV} to obtain
\begin{equation}
\label{franchetta}
\GDCH^*_B(S_b \times S_b) = \bigl\langle p_i^*(h_j)\bigr\rangle \oplus \QQ[\Delta_{S_b}]\ .
\end{equation}

This indeed uses the relations computed in \cite{BV}, and recalled in Section \ref{s:bv} for the reader's convenience, which express $\Delta^2_{S_b}$ and $\Delta_{S_b} h$ in terms of $h$ and the Beauville-Voisin class $o_{S_b}$. Third, the RHS of \eqref{franchetta} injects into 
$$
H^0(S)\otimes H^4(S) \oplus \CH^1(S)_{}\otimes \CH^1(S)\oplus H^4(S)\otimes H^0(S)
 \oplus {\mathbb Q}[\Delta_{S_b}] \subset H^4(S_b \times S_b),
$$
thus proving the Franchetta property for this family. This injection holds because the diagonal is linearly independent from "decomposable cycles" in $H^4(S_b \times S_b)$. More specifically, if the diagonal were decomposable, the homomorphism induced on $H^{2,0}(S_b)$ by the diagonal would be trivial, which can not be true for a K3 surface. 

We observe that the Cayley-Oguiso automorphism $g_b$ is defined on the whole $18$-dimensional family. In fact, the surfaces $S_i$ (and the maps $\phi_i$) in Section 3 exist relatively over $B$, and so the description
  \[ g=\phi_0^{-1} \circ \phi_3\ \] 
  makes sense family-wise.  Hence there exists $G \in \hbox{Aut}({\mathcal S}/B)$ such that $G_{|S_b}=g_b$. Now let us look at the graph $\Gamma_g \in \CH^2({\mathcal S} \times _B {\mathcal S})$. Recall that the very general $S_b$ has Picard number 2.
The fact that $g_b$ is anti-symplectic thus translates into the fact that the class in cohomology of the graph can be fiberwise decomposed as follows
\begin{equation}
\label{diagonal}
\Gamma_{g_b}= - \Delta_{S_b}+ \sum_{j=1}^2 a_j D_1^j \times D_2^j \ \ \hbox{in}\  H^4(S_b \times S_b), \qquad a_j \in {\mathbb Q},\ \ \hbox{$b\in B$\ very\ general.}
\end{equation}
(Here $D_1^1,D_1^2$ denotes a basis of $NS(S_b)_{\QQ}$, and $D_2^1,D_2^2$ is a dual basis in $NS(S_b)_{\QQ}$.)

Notice that the identity in \eqref{diagonal} involves only generically defined cycles, i.e., each class is the restriction of a global class defined on the total space of the family. As a consequence of the Franchetta property established above, the equality \eqref{diagonal} holds true in $\CH^2(S_b \times S_b)$, for the very general $b\in B$. Moreover, Voisin's spread result  \cite[Lemma 3.1]{Vo} then implies that \eqref{diagonal} is true in $\CH^2(S_b \times S_b)$ for every $b \in B$. Applying the equality \eqref{diagonal} to $\CH^2_{hom}(S_b)$, it follows that $(g_b)^*=-id$ on $\CH^2_{hom}(S_b)$ for every $b \in B$. 
This proves the Theorem.
 \end{proof}
 
 \begin{remark}
Recall that $S$ is embedded in ${\mathbb P}^3$ via a polarization $D_0$. As proved in \cite{F+}, the automorphism $g$ maps $D_0$ to another polarization, which is called $D_3$ in loc. cit. In other words, the class $h$ is preserved by $g^*$ and, accordingly, the class $h^2$, as well as the Beauville-Voisin cycle $o_{S}$ in $CH^2(S)$ are preserved by $g^\ast$. This completely describes the action of $g^*$ on $\CH^2(S)$: the splitting in eigenspaces of $g^*$ agrees with the Beauville--Voisin splitting, i.e.
\[ \CH^2(S)^- = \CH^2_{hom}(S) \ ,\ \  \CH^2(S)^+={\mathbb Q}[o_{S}]\ .\]
  \end{remark}  
 
 \subsection{An alternative argument}
 \label{s:alter}
 
 Here we present an alternative proof of Theorem \ref{main} using Hilbert schemes and the Beauville involution.
 
 \begin{proof}(Alternative proof of Theorem \ref{main}) Let $S^{[2]}$ denote the Hilbert scheme of $S$, and let $g_{ind}\in\aut(S^{[2]})$ be the automorphism induced by $g\in\aut(S)$, where $g$ is the generator desscribed in Section \ref{s: auto}. The Hilbert scheme $S^{[2]}$ has a
 {\em multiplicative Chow--K\"unneth decomposition\/}, in the sense of Shen--Vial \cite{SV}; this is established in \cite[Chapter 13]{SV} (cf. also \cite{V6} and \cite{O+} for the case of $S^{[n]}$, $n$ arbitrary).
 This implies in particular that there is a splitting
   \begin{equation}\label{split} \CH^4(S^{[2]})= \CH^4_{(0)}(S^{[2]})\oplus  \CH^4_{(2)}(S^{[2]}) \oplus \CH^4_{(4)}(S^{[2]})\ , \end{equation}  
   where the pieces $\CH^4_{(i)}(S^{[2]})$ have a certain explicit description (cf. \cite[Part 2]{SV}, or alternatively \cite[Theorem 3.5]{LV} for a summary).
   
 There is a natural correspondence-induced surjection
    \begin{equation}\label{iso} \Gamma_\ast\colon  \CH^4_{(2)}(S^{[2]}) \ \xrightarrow{\cong}\ \CH^2_{hom}(S)\ .\end{equation}
  
This correspondence is simply the incidence relation inside $S^{[2]}\times S$ (i.e. pairs $x\in S^{[2]}$, $z\in S$ such that $Supp(x)=z$). The map above is surjective, because $\CH^2_{(2)}(S^{[2]})$ is generated by expressions 
  \[ \langle [x,o_S] - [y,o_S], x,y\in S\rangle\ ,\] 
  where $o_S$ is the Beauville--Voisin distinguished $0$-cycle on $S$, and $[x,y]\in S^{[2]}$ stands for the unordered pair of points $x,y\in S$. Hence, given any degree zero $0$-cycle in S, say
$
a= x_1+\ldots + x_r -y_1-\ldots-y_r,
$
the corresponding element
$
\sum_i \langle [x_i,o_S] - [y_i,o_S]\rangle 
$
maps to $a$, i.e., there is a surjection as claimed.

The action of an automorphism on $\CH^4$ respects the splitting \eqref{split} (e.g., one way to see that the piece $\CH^4_{(2)}(S^{[2]})$ is preserved is by invoking the equality
    $\CH^4_{(2)}(S^{[2]})= S_1 \CH^4(S^{[2]})\cap \CH^4_{hom}(S^{[2]})$, where $S_1\CH^4_{(2)}(S^{[2]})$ denotes the subgroup generated by points with rational equivalence orbit of dimension at least 1, cf. \cite{Vois}). It is readily checked that the isomorphism $\Gamma_\ast$ is compatible with automorphisms, in the sense that there is a commutative diagram
    \[ \begin{array}[c]{ccc}
                   \CH^4_{(2)}(S^{[2]}) & \xrightarrow{\Gamma_\ast}& \CH^2_{hom}(S)    \\
                   &&\\
                  \ \  \downarrow{\scriptstyle (g_{ind})_\ast} &&    \ \  \downarrow{\scriptstyle (g_{})_\ast}   \\
                  &&\\
                     \CH^4_{(2)}(S^{[2]}) & \xrightarrow{\Gamma_\ast}& \CH^2_{hom}(S) \ .\\
                     \end{array}\]
                     
                    Because of the isomorphism \eqref{iso}, it suffices to prove that
                    \begin{equation}\label{suffices} (g_{ind})_\ast=-\ide\colon\ \                      
                      \CH^4_{(2)}(S^{[2]}) \ \to\ \CH^4_{(2)}(S^{[2]}) \ .\end{equation}
                      Thanks to the relation of Proposition \ref{og}, we have
                      \[  (g_{ind})_\ast=     (\iota_2)_\ast (\iota_1)_\ast (\iota_0)_\ast\colon\ \ \CH^4(S^{[2]})\ \to\ \CH^4(S^{[2]})   \ ,\]
                      where the $\iota_j$ are Beauville involutions of $S^{[2]}$ with respect to isomorphic projective embeddings of $S$.
                      But it is known that
                      \[ \iota_\ast=-\ide\colon\ \        \CH^4_{(2)}(S^{[2]}) \ \to\ \CH^4_{(2)}(S^{[2]}) \ \]         
                      for the Beauville involution $\iota$ of the Hilbert square of any quartic K3 surface $S$ \cite{Lat0}, \cite{Lat1}. This implies \eqref{suffices}, and so ends the proof.
                        \end{proof}

 \section{Some more results}\label{otherfamilies}
 
 In this section, we prove the conjecture for two more 18-dimensional families of K3 surfaces with an infinite-order automorphism.
 
 \begin{theorem}\label{main2} Let $S_{m,a}$ be a quartic K3 surface as in Theorem \ref{lee}, and let $g$ be a generator of $\aut(S)\cong\ZZ$. Assume $(m,a)=(2,2)$ or $(m,a)=(2,3)$.
 Then Conjecture \ref{conj} is true for $S$, i.e. the action of $g$ on $\CH^2_{hom}(S)$ is the identity for $(m,a)=(2,2)$, and minus the identity for $(m,a)=(2,3)$.
     \end{theorem}

     \begin{proof} This is the combination of the following two propositions:
     
     \begin{proposition}\label{prop1} Let $\Ss\to B$ be a family of surfaces, with the following properties:
     
     \begin{enumerate}
     
     \item each $S_b$, $b\in B$ is a smooth projective surface of degree $d$ in $\PP^3$;
     
     \item for each $S_b$, $b\in B$, there exists a curve $C_b\subset\PP^3$ such that (a) the linear system ${\mathcal L}:=\PP H^0(\PP^3,\OO_{\PP^3}(d)\otimes {\mathcal I}_{C_b})$ is non-empty, and (b) the general surface in ${\mathcal L}$ is part of the family $\Ss\to B$.
     \end{enumerate}
     
     Then
       \[ \GDCH^2_{\mathcal L}(S_b\times S_b) \ \subset\ \QQ[\Delta_{S_b}] + \bigl\langle (p_i)^\ast \GDCH^\ast_{\mathcal L}(S_b)\bigr\rangle\ ,\]
       for any $b\in B$.
      \end{proposition}
      
      \begin{proposition}\label{prop2} Let $S_b$ be a surface $S_{m,a}$ as in Theorem \ref{main2}. Then there exists a curve $C_b\subset S_b$ verifying Proposition \ref{prop1}.
      \end{proposition}
    
    It is readily seen that these two propositions together prove theorem \ref{main2}. Indeed, for any $b\in B$ one can find a curve $C_b$ and a linear system $\mathcal L$ as in Proposition \ref{prop1}. Let us assume $S=S_b$ is very general in the linear system ${\mathcal L}$ (and so has minimal Picard number $\rho=2$). As in the proof of Theorem \ref{main}, we now consider the correspondence  
      \[ \Psi_b:=
\Gamma_{g_b} \pm \Delta_{S_b}+ \sum_{j=1}^\rho a_j D_1^j \times D_2^j \ \ \ \in \CH^2(S_b\times S_b)\ ,
\]
where the $\pm$ sign is a $+$ sign for $(m,a)=(2,2)$ and a $-$ sign for $(m,a)=(2,3)$, and the $a_j\in\QQ$ are such that $\Psi_b$ is homologically trivial.
(Here $D_1^1,\ldots,D_1^\rho$ denotes a basis of $NS(S_b)_{\QQ}$, and $D_2^1,\ldots,D_2^\rho$ is a dual basis in $NS(S_b)_{\QQ}$.)

The correspondence $\Psi_b$ is generically defined, and so Proposition \ref{prop1} applies. Since the correspondence $\Psi_b$ is homologically trivial, there is no contribution of the diagonal (indeed, the action of $\Psi_b$ on $H^{2,0}(S_b)$ is zero), and so one has that $\Psi_b$ is {\em decomposable\/}, i.e.
  \[ \Psi_b\ \in\  \bigl\langle (p_i)^\ast \GDCH^\ast_{\mathcal L}(S_b)\bigr\rangle\ ,\]     
  for the very general $b\in {\mathcal L}$. The spread lemma \cite[Lemma 3.1]{Vo} then implies that $\Psi_b$ is decomposable for any $b\in {\mathcal L}$ such that $S_b$ is smooth. Looking at the action of $\Psi_b$ on $\CH^2(S_b)$ (and remembering that decomposable correspondences act as zero on $\CH^2_{hom}(S_b)$), one finds that $g$ acts in the expected way on zero-cycles, proving Theorem \ref{main2}.
  
 Let us now prove Proposition \ref{prop1}. This is a variant on the ``stratified projective bundle argument'' from \cite{FLV} and \cite{FLV3}. Indeed, let $\Ss\to{\mathcal L}$ be the universal family. By assumption, the base locus of the linear system ${\mathcal L}$ consists of $C_b=C_b^0$ plus (possibly) some other curves $C_b^1,\ldots,C_b^r$.
  The fiber product
   \[ \Ss\times_{\mathcal L}\Ss\ \to\ \PP^3\times\PP^3 \]
   has the structure of a stratified projective bundle, with strata given by
   \[  C_b^i\times\PP^3\ ,\ \PP^3\times C_b^i\ ,\    \Delta_{\PP^3} \ ,\]
   plus the intersections of these loci. Applying \cite[Proposition 5.2]{FLV} and restricting to codimension 2, we find that
     \[ \begin{split} \GDCH^2_{\mathcal L}(S_b\times S_b) = &\bigl\langle (p_i)^\ast(h)\bigr\rangle + \QQ[\Delta_{S_b}]  +\sum_{i,j=0}^r\QQ[C^i_b\times C^j_b] \\
     &+ \sum_{i=0}^r\ima\bigl( \GDCH_2(C^i_b\times S_b)\to \CH^2(S_b\times S_b)\bigr)\\ &+  \sum_{i=0}^r\ima\bigl( \GDCH_2(S_b\times C^i_b)\to \CH^2(S_b\times S_b)\bigr)\ .\\
     \end{split}  \]
     Since $H^1(S_b)=0$, the last two terms are decomposable (cf. Lemma \ref{q} below). This proves the proposition.
     
     \begin{lemma}\label{q} Let $Y_1,Y_2$ be two smooth projective varieties, and assume that $H^1(Y_1,\QQ)=0$. Then
       \[ \CH^1(Y_1\times Y_2)=  (p_1)^\ast \CH^1(Y_1) \oplus (p_2)^\ast \CH^1(Y_2)\ ,\]
       where $p_i$, $i=1,2$ are the two projections.
     \end{lemma}
     
    \begin{proof}(of the lemma) This is surely well-known. The K\"unneth formula implies
      \[ \begin{split} H^2(Y_1\times Y_2,\QQ)&= (p_1)^\ast H^2(Y_1,\QQ)\oplus (p_2)^\ast H^2(Y_2,\QQ)  \ ,\\
      \  H^2(Y_1\times Y_2,\OO)&= (p_1)^\ast H^2(Y_1,\OO)\oplus (p_2)^\ast H^2(Y_2,\OO)\ ,\\
      \end{split}\]
      and hence
      \begin{equation}\label{ns}NS(Y_1\times Y_2)_{\QQ}= (p_1)^\ast NS(Y_1)_\QQ\oplus (p_2)^\ast NS(Y_2)_\QQ\ .\end{equation}
      Likewise, one has
      \[ H^1(Y_1\times Y_2,\QQ)=  (p_2)^\ast H^1(Y_2,\QQ)     \ ,\]
      and so
               \begin{equation}\label{pic0} \pic^0(Y_1\times Y_2)_{\QQ}=  (p_2)^\ast \pic^0(Y_2)_\QQ\ .\end{equation}
      Combining equalities \eqref{ns} and \eqref{pic0}, we find that $\CH^1(Y_1\times Y_2)=\pic(Y_1\times Y_2)_{\QQ}$ is decomposable, proving the lemma.
      \end{proof}
      
Let us now prove Proposition \ref{prop2}. As proved in \cite[\S 5.1]{Lee}, the ample cone of $S_{(2,a)}$ is given by curves $C=x+y\theta$, where $\theta$ satisfies $\theta^2=1+a\theta$, the integer $x$ is positive and $y$ is a positive integer such that $C^2=4(x^2+axy-y^2)>0$. Moreover, any such divisor is very ample and, accordingly, effective. In order to define the pair $(C_b, S_b)$ that satisfies the assumptions of the proposition, we need a polarization $H_b$ such that $H^2_b=4$ and a divisor $\Gamma_b$. The former gives an embedding of $S_{(2,a)}$ to a smooth quartic surface $S_b$ in ${\mathbb P}^3$; the latter maps $\Gamma_b$ to $C_b$.

If $a=2$, set $H_b=1+2\theta$ and $\Gamma_b=1+\theta$. For $a=3$ choose $1+3\theta$ and $2+7\theta$, respectively. It is easy to check that in both cases $H_b^2=4$. This proves that  the assumption of Proposition \ref{prop1} (1) is verified for $d=4$. As for the assumption (2) in Proposition \ref{prop1}, let us consider the short exact sequence
\begin{equation}
0 \to {\mathcal O}_{{\mathbb P}^3}(4) \otimes {\mathcal I}_{C_b} \to {\mathcal O}_{{\mathbb P}^3}(4) \to {\mathcal O}_{C_b}(4) \to 0.
\end{equation}  

The long exact sequence in cohomology gives
\begin{equation}
h^0({\mathbb P}^3,{\mathcal O}_{{\mathbb P}^3}(4) \otimes {\mathcal I}_{C_b}) = 35 - h^0(C_b, {\mathcal O}_{C_b}(4) ) + h^1({\mathbb P}^3, {\mathcal O}_{{\mathbb P}^3}(4) \otimes {\mathcal I}_{C_b}).
\end{equation}                  

If we prove that our choice of $\Gamma_b$ and $H_b$ implies the inequality
\begin{equation}
h^0(C_b, {\mathcal O}_{C_b}(4) ) \leq 33 + h^1({\mathbb P}^3, {\mathcal O}_{{\mathbb P}^3}(4) \otimes {\mathcal I}_{C_b}),
\end{equation}
we will have $h^0({\mathbb P}^3,{\mathcal O}_{{\mathbb P}^3}(4) \otimes {\mathcal I}_{C_b}) \geq 2$ and so assumption (2) of Proposition \ref{prop1} is satisfied. 

By direct computation, we have 
\[ \deg\left(K_{C_b} \otimes {\mathcal O}_{C_b}(-4)\right)=\Gamma_b^2-4\Gamma_bH_b<0 \] for $a=2,3$. Hence the Riemann-Roch theorem for $C_b$ implies that
$$
h^0(C_b, {\mathcal O}_{C_b}(4))=4\Gamma_b H_b - \frac{\Gamma_b^2}{2}
$$
for $a=2,3$.

Our choice of $\Gamma_b$ and $H_b$ for $a=2,3$ shows that $4\Gamma_b H_b - \frac{\Gamma_b^2}{2} \leq 33$ (direct computation), and so we have ascertained that Proposition \ref{prop2} can be applied to this set-up, concluding the proof.

\end{proof}
     
     \begin{remark} The argument of this section also applies to the Cayley--Oguiso surfaces, giving a third proof of Theorem \ref{main}. Indeed,
     the Cayley--Oguiso determinantal quartics are characterized as those quartics in $\PP^3$ containing a sextic curve of genus 3.
     
     Note also that the alternative proof given in Subsection \ref{s:alter} does {\em not\/} apply to the surfaces of Theorem \ref{main2}. This is because for these surfaces the relation with the Beauville involution (Proposition \ref{lee2}) is different in nature from the relation with the Beauville involution for the Cayley--Oguiso surfaces (Proposition \ref{og}).
      \end{remark}

 \vskip1.5cm
 
 \vskip0.5cm
\begin{nonumberingcon} The authors declare that there are no conflicts of interest related to this work.
\end{nonumberingcon}

\vskip0.5cm

 \begin{nonumberingt} 
Thanks to Pasticceria Vabres in Palermo for their wonderful arancine. Thanks to the referee for pertinent and helpful comments.
\end{nonumberingt}

\end{document}